\documentclass[12pt,a4paper]{amsart}

\usepackage{amssymb,amsmath,amsthm}
\usepackage{color}
\usepackage{rotating}
\usepackage{fullpage}
\usepackage{comment}
\usepackage{hyperref}
\usepackage{soul}  
\usepackage[all]{xy}
\DeclareMathOperator{\inte}{int}

\newtheorem{theorem}{Theorem}
\newtheorem{lemma}[theorem]{Lemma}

\theoremstyle{definition}
\newtheorem{definition}{Definition}
\newtheorem{example}{Example}

\title{Periodic orbits 
in R\"ossler system}

\author{Anna Gierzkiewicz}
\email{anna.gierzkiewicz@urk.edu.pl}
\address{Department of Applied Mathematics, University of Agriculture in Krak\'ow,
ul. Balicka 253c,
30--198 Krak\'ow, Poland
}

\author{Piotr Zgliczy\'nski}
\email{umzglicz@cyf-kr.edu.pl}
\address{Institute of Computer Science, Jagiellonian University,
ul. \L ojasiewicza 6, 30-348 Krak\'ow, Poland
}

\begin{document}

\begin{abstract}

We prove the existence of $n$-periodic orbits for almost all $n\in\mathbb{N}$ in the R\"ossler system with attracting periodic orbit, for two sets of parameters. The proofs are computer-assisted.


\end{abstract}

\maketitle


\section{Introduction}
R\"ossler in \cite{Rossler76} studied the system
\begin{equation}\label{eq:rossler}
\begin{cases}
x'=-y-z,
\\
y'=b y+x,
\\
z'=z (x-a)+b
\end{cases}
\end{equation}
with $a=5.7$, $b=0.2$ as an example of a low-dimensional polynomial system with a single nonlinear term, which admits chaotic behaviour. The existence of symbolic dynamics in it was proven with computer assistance in \cite{ZRossler}. Studying R\"ossler system's dynamics and its periodic orbits for varying parameters is an active field of research (see \cite{Glen, Kri, TP, Valls, Pil, Galias} and references therein).

Numerical simulations of the R\"ossler system show that for wide range of parameters for suitable Poincar\'e map there exists an attracting set, which is essentially one dimensional. In such a situation one expects that the forcing relations between periods described in the Sharkovskii's theorem should be applicable.
  Let us recall the Sharkovskii's theorem in its original form in dimension one (see \cite{ShU,Stefan,Block}):

\begin{theorem}[Sharkovskii]\label{th:shar}
Define an ordering `$\triangleleft$' of natural numbers:
\[
\begin{array}{r@{\ \triangleleft\ }c@{\ \triangleleft\ }c@{\ \triangleleft\ }c@{\ \triangleleft\ }l}
3\triangleleft 5 \triangleleft 7 \triangleleft 9 \triangleleft\ \dots & 2\cdot 3 & 2 \cdot 5 & 2 \cdot 7 & \dots
\\
 \dots & 2\cdot 3^2 & 2 \cdot 5^2 & 2 \cdot 7^2 & \dots
\\
 \dots & 2\cdot 3^3 & 2 \cdot 5^3 & 2 \cdot 7^3 & \dots
\\
 \multicolumn{5}{c}{\dots \dots \dots \dots \dots \dots \dots \dots \dots \dots }
\\
 \dots & 2^{k+1} & 2^k & 2^{k-1} & \ldots\ \triangleleft 2^2 \triangleleft 2 \triangleleft 1.
\\
\end{array}
\]
Let $f: I \to \mathbb{R}$ be a continuous map of an interval. If $f$ has an $m$-periodic point and $m\triangleleft n$, then $f$ also has an $n$-periodic point.
\end{theorem}

In \cite{PZszarI,PZszarII} the above theorem has been extended to multidimensional perturbations of 1-dimensional maps and in \cite{PZmulti} it is explained how it fits to the symbolic dynamics established  \cite{ZRossler} for system (\ref{eq:rossler}) with  $a=5.7$, $b=0.2$.

\begin{figure}[h]
	\includegraphics[height=7cm]{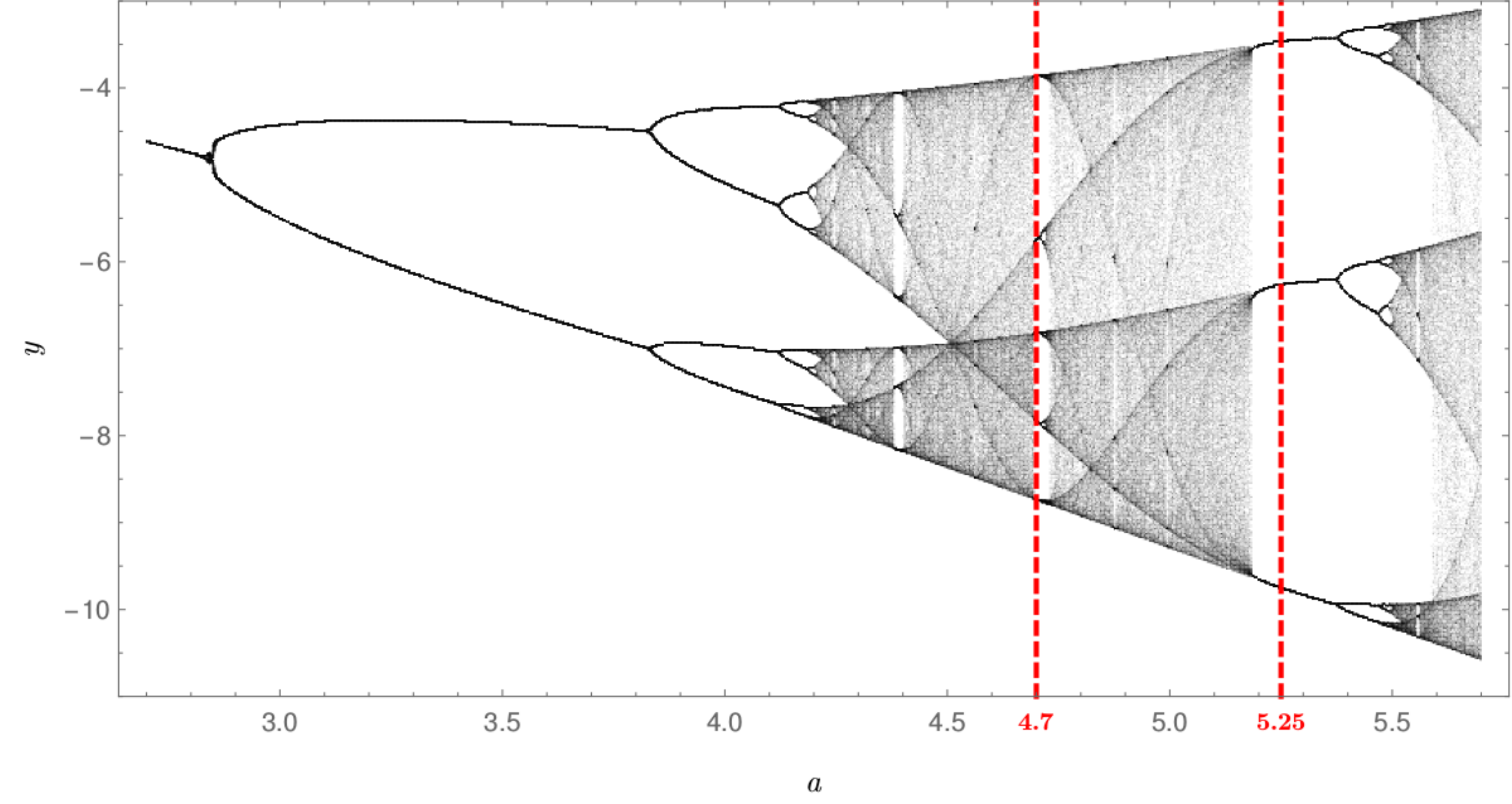}
	\caption{\label{fig:bif}Bifurcation diagram for the Poincar\'e map of the system \eqref{eq:rossler} with $b=0.2$. The attracting periodic orbits for $a=4.7$ and $a=5.25$ are expected ($3$-periodic for $a=5.25$ and $5$-periodic for $a=4.7$, respectively). }
\end{figure}

In our study of (\ref{eq:rossler}) we fix $b=0.2$. Consider the Poincar\'e map $P$ partly defined on the section $\Pi=\{x=0, y<0\}$. Studying the bifurcation diagram for $(a,y)$ (see Fig. \ref{fig:bif}) we see that  a $3$-periodic attracting orbit appears, for example, for $a=5.25$. The other interesting case is the $5$-periodic attracting point for $a=4.7$. These are two sets of parameter values treated in our work.

From numerical evidence we see that $P$ exhibits a strong contraction in the $z$ direction and the invariant sets containing the periodic orbits are almost 1-dimensional (see Figs. \ref{fig:intervalp123} and \ref{fig:intervalp12345}). Therefore we can assume in both cases considered in our work (and, most probably, also for a larger range of $a$s) that $P$ can be treated as a 2-dimensional perturbation of a 1-dimensional \textit{model map}.  Hence, as suggested by Theorem~\ref{th:shar} and its multidimensional version from \cite{PZszarI,PZszarII,PZmulti}, we can expect a large set of periodic points also for $P$.

In this paper we prove the existence of all periods for $a=5.25$ (with an attracting $3$-periodic orbit). In the second case, $a=4.7$ (with  an attracting orbit of period $5$), we show the existence of all periods with an exception of the period $3$, which agrees with the forcing relation between periods for interval maps described in Theorem~\ref{th:shar}. Moreover, we show that the period $3$ is indeed not realised in an attracting set  for $P$.

The methods used in the present work combine the topological tools, the covering relations, with rigorous numerics. Our proofs are inspired
by the approach to multidimensional perturbations of 1-dimensional maps from \cite{PZszarI,PZszarII,PZmulti}, but the strategy of constructing
suitable covering relations needed to obtain the desired infinite sets of periodic orbits is different. In the second case some orbits
with low periods has been established using the interval Newton method \cite{N}.

All the programs are written in C++ with the use of CAPD (Computer-Assisted Proofs in Dynamics) library \cite{capd,capd-article} for interval arithmetic, differentiation and integration.

\section{Notation}

In the paper, we consider the system \eqref{eq:rossler} with fixed $b=0.2$ and $a=5.25$ or  $a=4.7$. In both cases, we denote by $\Pi$ the half-plane $\{x=0, y<0\}$ with induced coordinates $(y,z)$, and $P$ is a Poincar\'e map on section $\Pi$, that is the map
\[
P(y,z) = \pi_{(y,z)}\left( \Phi_{T(y,z)}\left(x=0,y,z\right)\right)\text{,}
\]
where $\pi_{(y,z)}$ is the projection on the $(y,z)$ plane, $\Phi_t$ is the dynamical system induced by considered system and $T=T(y,z)$ is a return time, if well-defined. Note that for $y+z=0$ the vector field given by the right-hand-side of \eqref{eq:rossler} is not transversal to $\Pi$. In our area of interest, however, $z$ is sufficiently small to guarantee $\dot{x}>0$ on the section.

To simplify the notation, by `$n$-periodic orbit' or `point', we understand an orbit or point with basic period $n$ for map $P$.
Whenever we refer to a `$n$-periodic orbit of the system' we mean a periodic orbit of the system, which passes through an $n $-periodic orbit of $P$.

\section{Horizontal covering and periodic points}\label{sec:covering}

The idea of one-dimensional covering relation between intervals which Block \textit{et al.} \cite{Block} used to prove Sharkovskii's theorem (Th. \ref{th:shar}) is given by the following definition.
\begin{definition}
Assume that $I, J \subset \mathbb{R}$ are intervals and $f:I\to \mathbb{R}$ is continuous.

An interval $I$ \emph{$f$-covers} $J$ (denoted by $I \overset{f}{\longrightarrow} J $) if there exists a subinterval $K\subset I$ such that $f(K)=J$.
\end{definition}

The following easy theorem gives the periodic orbits in the proofs of Sharkovskii's theorem (see, for example, \cite{Block}).
\begin{theorem}
\label{th:1d-covering}
Assume that $f:\mathbb{R} \supset I \to \mathbb{R}$ is continuous and we have a sequence of intervals $I_j \subset I$ for $j=0,\dots,n-1$ such that
\begin{equation*}
  I_0  \overset{f}{\longrightarrow} I_1 \overset{f}{\longrightarrow} I_2 \overset{f}{\longrightarrow} \dots \overset{f}{\longrightarrow} I_{n-1} \overset{f}{\longrightarrow} I_0.
\end{equation*}
Then there exists $x \in I_0$, such that $f^j(x) \in I_j$ for $j=1,\dots,n-1$ and $f^n(x)=x$.
\end{theorem}

For higher-dimensional perturbations of $f$ a stronger notion of covering is needed to have an analogous result.
We use the notion of the horizontal covering from \cite{PZszarI,PZmulti} with small modifications. It is simplified, set in two-dimensional space, and the notion of `left-' and `right side' is slightly extended.

By $\mathcal{C}(r)$ we will denote the family of rectangles (two-dimensional cylinders) of the form $[a,b]\times [-r,r]\subset \mathbb{R}^2$ for $r>0$.

\begin{definition}
A \textit{two-dimensional h-set} (shortly: an h-set) is a rectangle $N = [a,b]\times [-r,r] \in \mathcal{C}(r)$ with the following elements distinguished:
\begin{itemize}
	\item its left edge \quad $L(N) = \{a\} \times [-r,r]$;
	\item its right edge \quad $R(N) = \{b\} \times [-r,r]$;
	\item its horizontal boundary \quad $H(N) = [a,b] \times \{-r,r\}$;
	\item its left side \quad $S_L(N) = (-\infty, a) \times (-\infty,\infty)$;
	\item its right side \quad $S_R(N) = (b,\infty) \times (-\infty,\infty)$.	
\end{itemize}
\end{definition}

We need these notions to define horizontal covering relation.

\begin{figure}[h]
	\includegraphics[height=7cm]{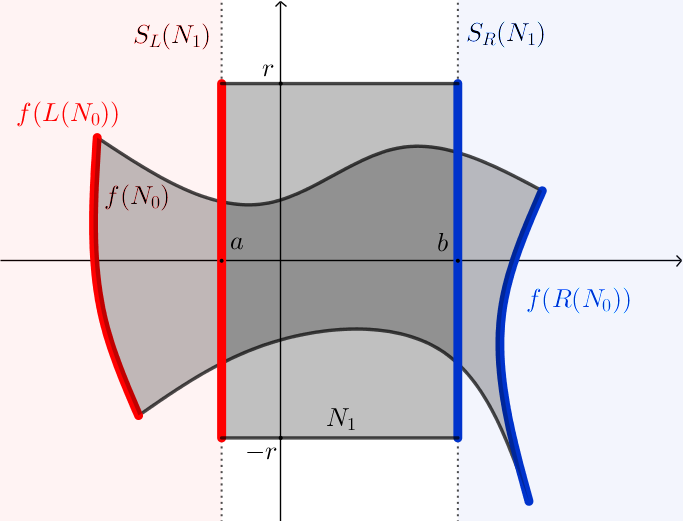}
	\caption{\label{fig:cover}Horizontal covering $N_0 \overset{f}{\Longrightarrow} N_1$}
\end{figure}

\begin{definition}
\label{def:cov}
Let $N_0$, $N_1 \in \mathcal{C}(r)$ be two h-sets and $f: \mathbb{R}^2 \to \mathbb{R}^2$.
We say that $N_0$ \emph{$f$-covers $N_1$ horizontally} and denote by $N_0 \overset{f}{\Longrightarrow} N_1$ if
\begin{equation}
f(N_0) \subset S_L(N_1)\cup N_1 \cup S_R(N_1) \setminus H(N_1)\text{,}  \label{eq:cov-miedzy}
\end{equation}
and one of the two conditions hold:
\begin{equation}
\begin{aligned}
&\text{either } & f(L(N_0))\subset S_L(N_1) &\text{\quad and \quad} f(R(N_0))\subset S_R(N_1)\text{,}
\\
&\text{or } & f(L(N_0))\subset S_R(N_1) &\text{\quad and \quad} f(R(N_0))\subset S_L(N_1).
\end{aligned}
\end{equation}
\end{definition}
See Fig. \ref{fig:cover} for the illustration of horizontal covering.
Let us emphasize that the above conditions can be easily checked with the use of computer via interval arithmetic and `$<$', `$>$' relations.

The following theorem might be seen as the generalization of Theorem~\ref{th:1d-covering}.

\begin{theorem}[\cite{PZszarI}]\label{th:periodic}
Suppose that we have a loop of $n$ horizontal $f$-coverings:
\[
N_0 \overset{f}{\Longrightarrow} N_1\overset{f}{\Longrightarrow} \dots \overset{f}{\Longrightarrow} N_{n-1} \overset{f}{\Longrightarrow} N_n = N_0\text{,}
\]
then there exists $x\in \inte N_0$ such that $f^n(x)=x$ and
\[
\text{for } i=0,\dots ,n-1: \qquad f^i(x)\in \inte N_i.
\]
\end{theorem}

The example of topological horseshoe discussed below shows how from a finite number of covering relations we can obtain periodic orbits of all periods.
\begin{example}[A topological horseshoe]
Let $N_0$, $N_1 \subset \mathcal{C}(r)$ be two disjoint h-sets. Suppose that a continuous map $f:\mathbb{R}^2 \to \mathbb{R}^2$ fulfils the horizontal covering relations (see Fig. \ref{fig:horseshoe})
\begin{equation}\label{eq:horseshoe}
\begin{array}{cc}
N_0 \overset{f}{\Longrightarrow}N_0\text{,} \quad
&N_0 \overset{f}{\Longrightarrow}N_1\text{,}
\\
N_1 \overset{f}{\Longrightarrow}N_0 \text{,} \quad
&N_1 \overset{f}{\Longrightarrow}N_1.
\end{array}
\end{equation}

\begin{figure}[h]
	\includegraphics[height=5cm]{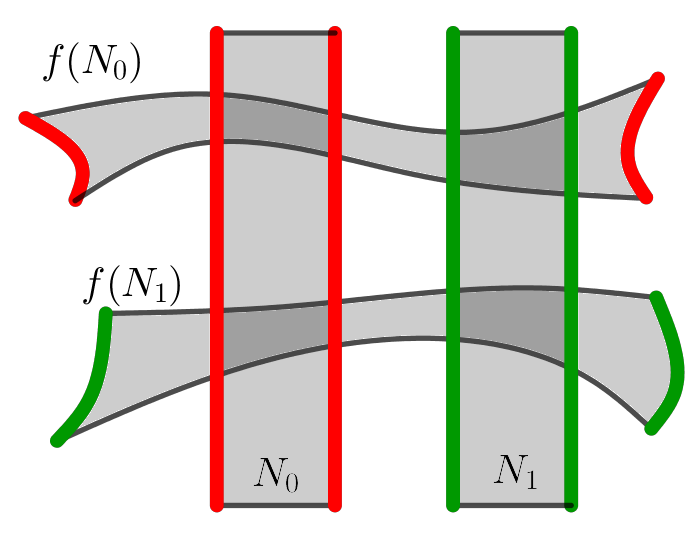}
	\caption{\label{fig:horseshoe}Topological horseshoe in $\mathbb{R}^2$: each $N_{0,1}$ covers itself and the other set. The vertical egdes of $N_0$ and $N_1$ are marked in red and green, respectively.}
\end{figure}
Such a map is called in literature a \emph{topological horseshoe} for $N_0$, $N_1$ \cite{Morse}.

Choose now any finite sequence of zeros and ones of any length $l$ : $(i_0, i_1, \dots, i_{l-1})$, $i_k\in\{0,1\}$. The conditions \eqref{eq:horseshoe} imply in particular the following chain of covering relations:
\[
N_{i_0} \overset{f}{\Longrightarrow}
N_{i_1} \overset{f}{\Longrightarrow}
\dots \overset{f}{\Longrightarrow}
N_{i_{l-1}} \overset{f}{\Longrightarrow}
N_{i_0}\text{,}
\]
and, from Theorem \ref{th:periodic}, we can deduce that there exists an $l$-periodic orbit for $f$, which additionally moves between the sets $N_0$, $N_1$ according to the pattern
$(i_0, i_1, \dots, i_{l-1})$. We can do it for any period $l$ and we are able to choose such a sequence of length $l$ to be sure that $l$ is the fundamental period of the orbit.

The same is also true for (bi-)infinite sequences of indicators $i_k$. In general, the maps admitting a topological horseshoe are semi-conjugate to the dynamical system generated on the space $\Sigma_2 = \{0,1\}^{\mathbb{Z}}$ by the shift map, called also the \textit{symbolic dynamics} \cite{Morse}.


\end{example}

\section{Orbits of all periods for $a=5.25$}

Consider now the system \eqref{eq:rossler} with $a=5.25$, $b=0.2$, that is
\begin{equation}\label{eq:rossler525}
\begin{cases}
x'=-y-z,
\\
y'=0.2\, y+x,
\\
z'=z (x-5.25)+0.2
\end{cases}
\end{equation}
 We expect from the bifurcation diagram (Fig. \ref{fig:bif}) that there exists a $3$-periodic orbit for the system \eqref{eq:rossler525}. The Lemma \ref{lem:3-per} below establishes this fact.

\begin{lemma}\label{lem:3-per}
The Poincar\'e map $P$ of the system \eqref{eq:rossler525} has a 3-periodic orbit $\mathcal{O}^3=\{p_1^3,p_2^3,p_3^3\}$ contained in the following rectangles in the $(y,z)$ coordinates on the section $\Pi$:
\begin{equation}\label{eq:3-per}
	\begin{aligned}
		p_1^3 \in & -3.4664152050_{12922}^{08744} \times 0.034631605476_{4013}^{51117}\text{,} \\	
		p_2^3 \in & -6.2640075332_{82922}^{74157} \times 0.032654358846_{02701}^{20798}\text{,} \\
		p_3^3 \in & -9.7488899180_{93569}^{88608} \times 0.030752873380_{62635}^{70747}\text{,}
		\end{aligned}
\end{equation}
\end{lemma}

\begin{figure}[h]
	\includegraphics[height=7cm]{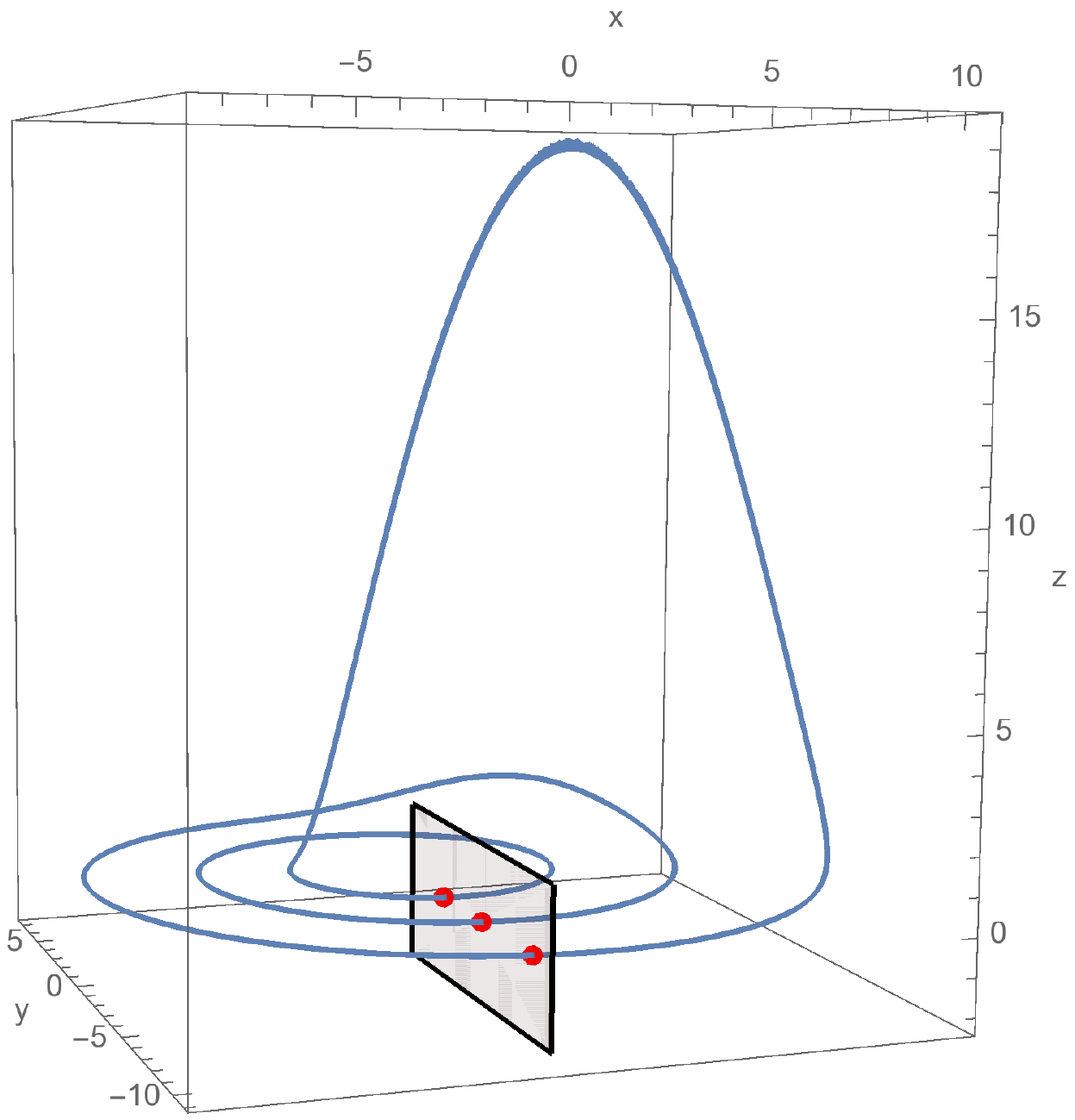}
	\caption{\label{fig:attr3} The attracting $3$-periodic orbit for the system \eqref{eq:rossler525}.}
\end{figure}

\begin{proof}
Computer-assisted by interval Newton method applied to $P$: \cite{proof}, Case 1 of the program \texttt{01-Roessler\_a525.cpp}. See also the outline in Appendix, Subsec.\ \ref{ss:prog_periodic}.
\end{proof}

\begin{theorem}\label{th:r3}
The R\"ossler system \eqref{eq:rossler525} has $n$-periodic orbits for any $n\in \mathbb{N}$.
\end{theorem}


Before we present a formal proof we explain the heuristics behind our construction.

\begin{figure}[h]
	\includegraphics[width=0.45\textwidth]{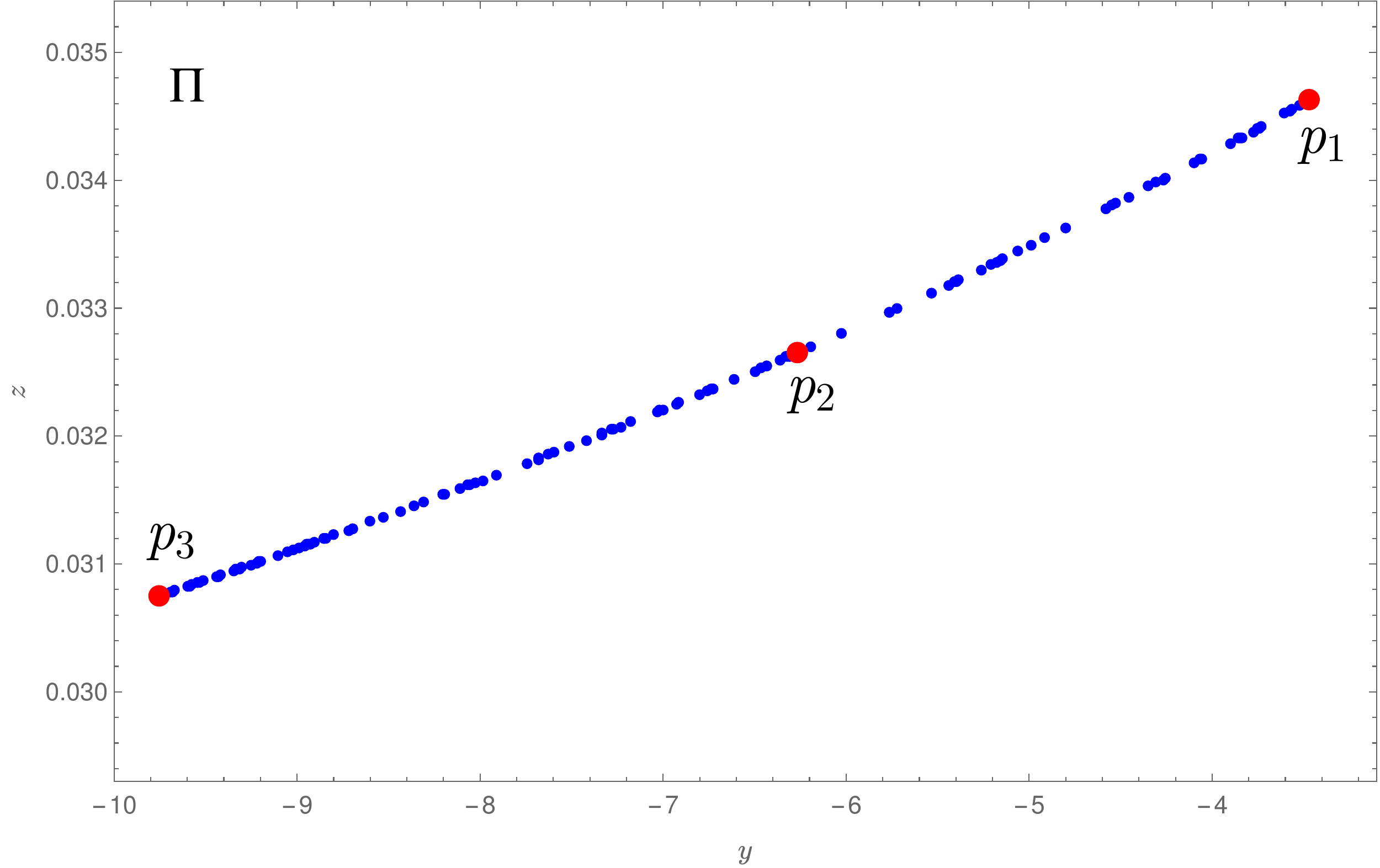}	
	\quad	
	\includegraphics[width=0.45\textwidth]{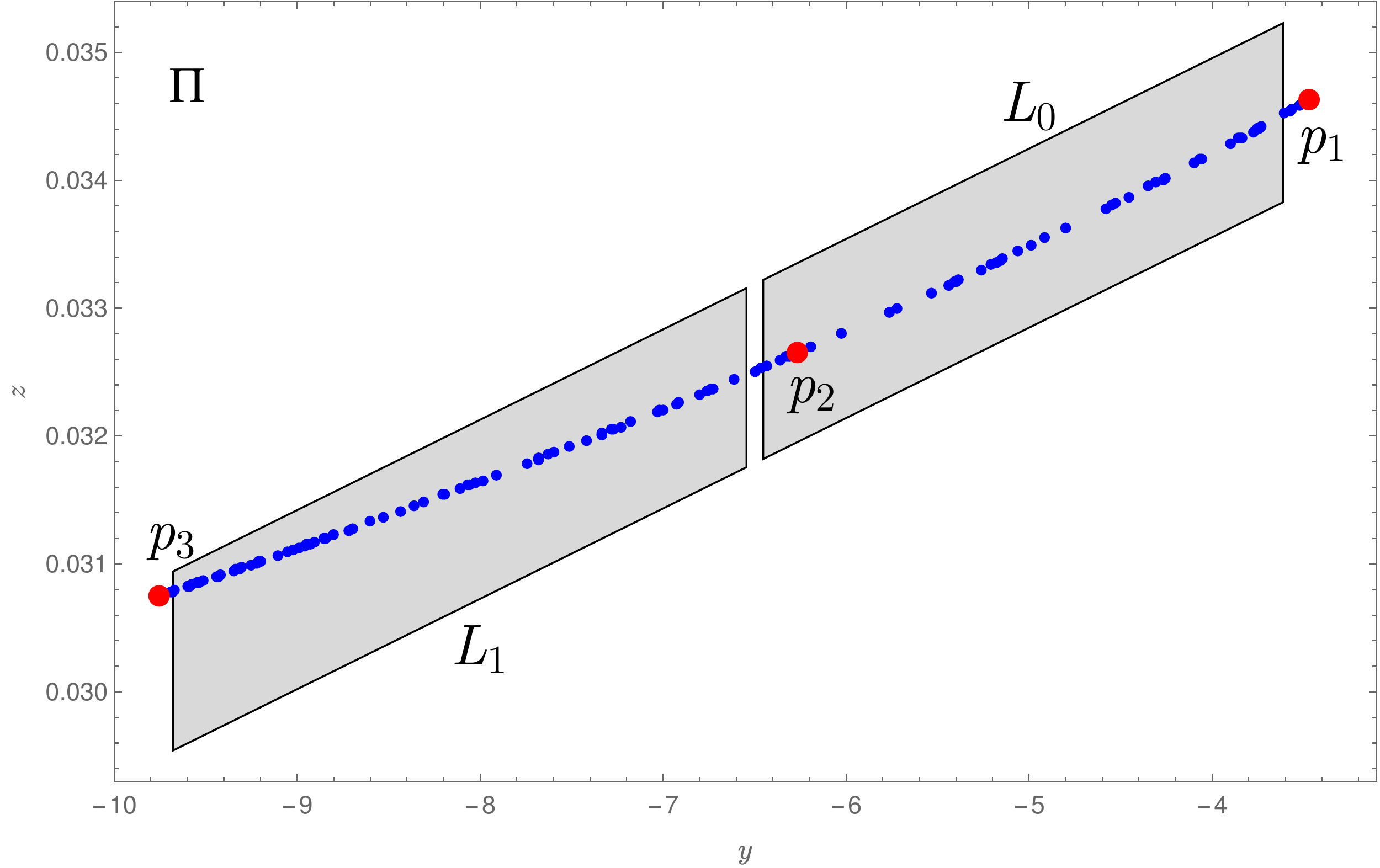}	
	\caption{\label{fig:intervalp123}\underline{To the left:} the fragment of $\Pi$ section containing the $3$-periodic orbit $\mathcal{O}^3$ (red) and some orbits attracted by it. The attracted orbits sketch the shape of the hypothetical $1$-dim manifold $\mathcal{M}$.
	\newline
	\label{fig:L_vs_p3}\underline{To the right:} the location of the sets $L_0$, $L_1$ relative to the orbit $\mathcal{O}^3$.}
\end{figure}

We follow the idea of the proof of \cite[Th. 2.11]{PZmulti}. Assume that there exists a 1-dimensional manifold $\mathcal{M}\subset \Pi$ with boundary 
 (or simply a homeomorphic image of a closed interval), such that $\mathcal{O}^3\subset \mathcal{M}$ and apparently $P(\mathcal{M})$ is contained in a small neighbourhood of $\mathcal{M}$ (see Fig. \ref{fig:intervalp123}).

\begin{figure}[h]
	\includegraphics[width=0.5\textwidth]{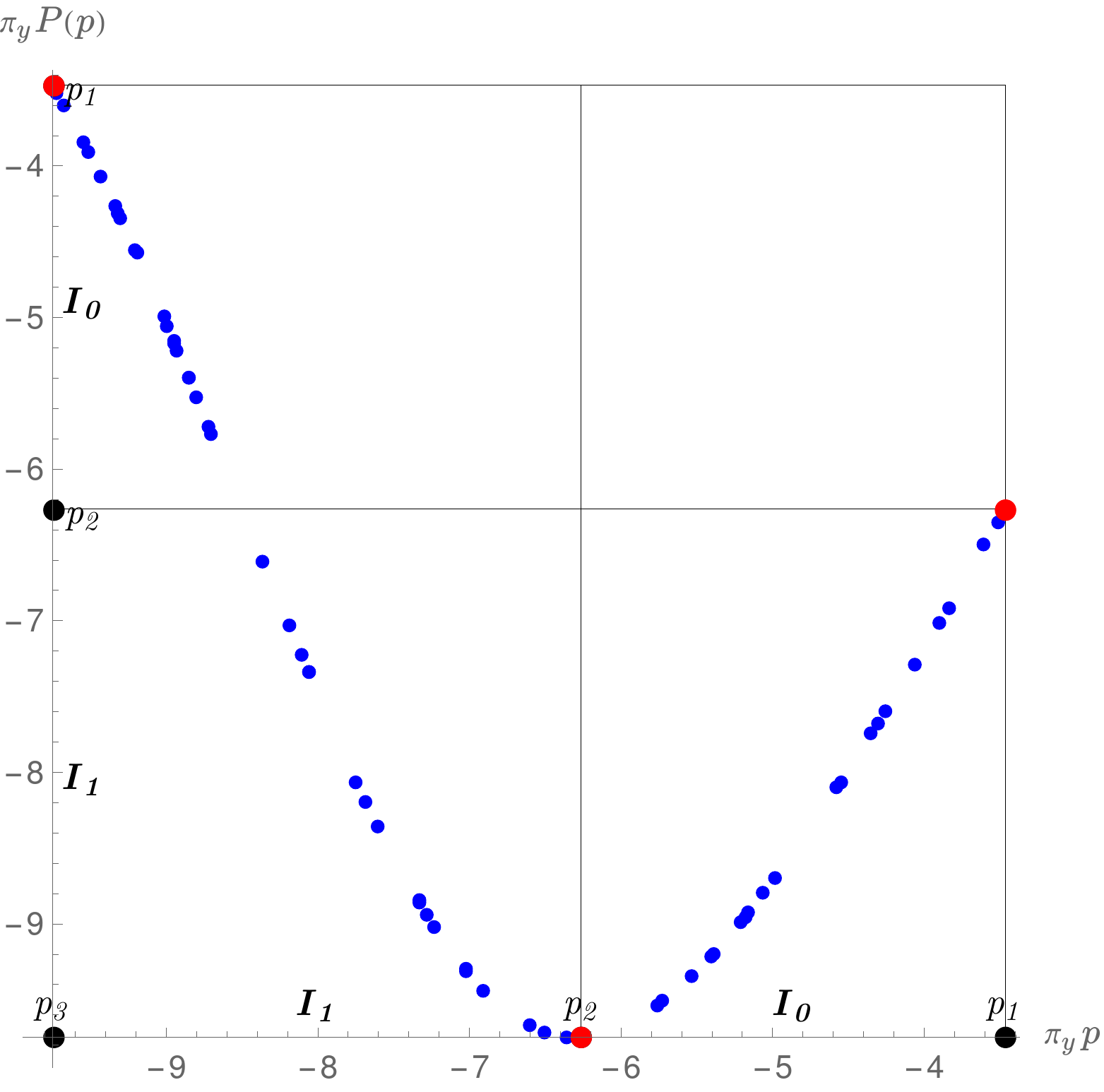}
	\caption{\label{fig:py_vs_Ppy}The images of the points near $\mathcal{M}$ through $P$ sketch the plot of the model map $\mathcal{P}$.}
\end{figure}

Considering Fig. \ref{fig:intervalp123}, we see that it is reasonable to parameterize $\mathcal{M}$ by the $y$ coordinate. We can now try to make a numerical `plot' of $\mathcal{M}$'s self  map  $\mathcal{P}$ (a model map for $P$), using the $y$ coordinates of the attracted orbits:
\[
y \mapsto \pi_y P(y,z)=:\mathcal{P}(y) \text{,\quad where $(y,z)$ apparently belongs to $\mathcal{M}$,}
\]
as on Fig.\ \ref{fig:py_vs_Ppy}.

Let us denote the `segment' $[p_1^3,p_2^3]$ by $I_0$ and $[p_2^3,p_3^3]$ by $I_1$, as on Fig.\ \ref{fig:py_vs_Ppy}. Note that $\mathcal{M}$'s segments would fulfil the one-dimensional covering relations:
\[
I_0 \overset{\mathcal{P}}{\longrightarrow} I_1 \overset{\mathcal{P}}{\longrightarrow} I_1 \overset{\mathcal{P}}{\longrightarrow}I_0 .
\]
Basing on this observation, we find two-dimensional sets $L_0$, $L_1$ on $\Pi$, lying close to the segments $I_0$, $I_1$, which most probably fulfil the similar horizontal covering relations. Their location relative to the orbit $\mathcal{O}^3$ is depicted on Fig. \ref{fig:L_vs_p3}, to the right.
In the proof below we show that these horizontal covering relations indeed occur.

\begin{proof}~

Let $M = \begin{bmatrix}
		-1 & 0.000706767 \\-0.000706767 & -1
	\end{bmatrix}$ and (with some abuse of notation) $p_2^3=\begin{bmatrix}
	-6.264007533274157 \\ 0.03265435884602701
	\end{bmatrix}$. Define an affine map $C$ on $\mathbb{R}^2$ by
	\[
	C(x)=Mx+p_2^3.
	\]
The matrix $M$ is chosen to place the images of horizontal h-sets through $C$ approximately along $\mathcal{M}$.
Consider now two h-sets $N_0,N_1\subset \mathbb{R}^2$:

\begin{equation}
\begin{aligned}
	N_0 &= [	-1.23094 \pm 1.41278]\times [\pm 7\cdot 10^{-4}]
	\text{,}
	\\
	N_1 &= [	1.84699 	\pm 1.55949]\times [\pm 7\cdot 10^{-4}].
\end{aligned}
\end{equation}
Denote their images through $C$ by
\[
L_0 = C(N_0)\text{, }\qquad L_1 = C(N_1)\text{, }
\]
and consider these parallelograms as sets on the section $\Pi$ (see Fig. \ref{fig:3-per}).

\begin{figure}[h]
	\includegraphics[width=0.45\textwidth]{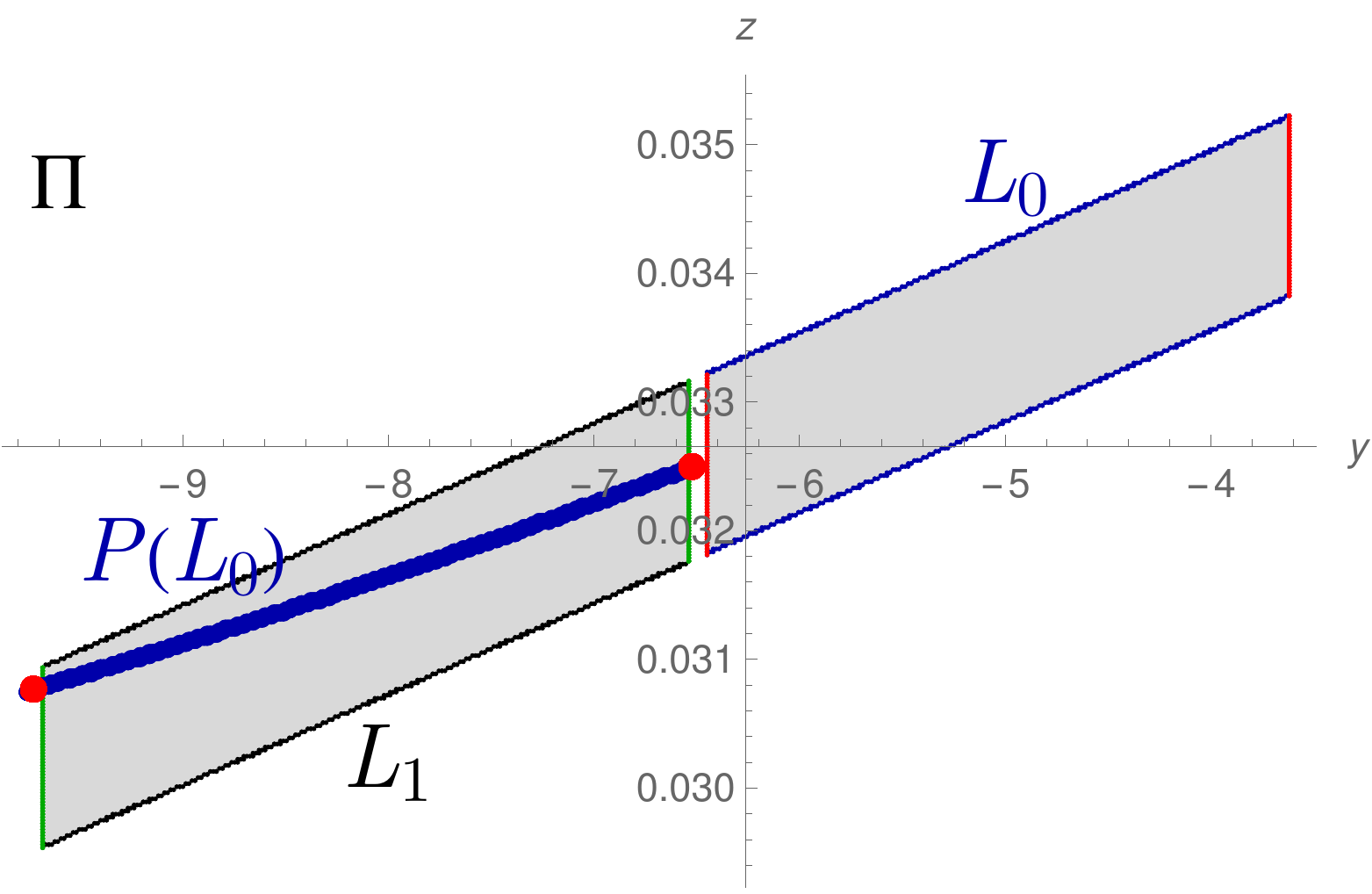}	
	\quad
	\includegraphics[width=0.45\textwidth]{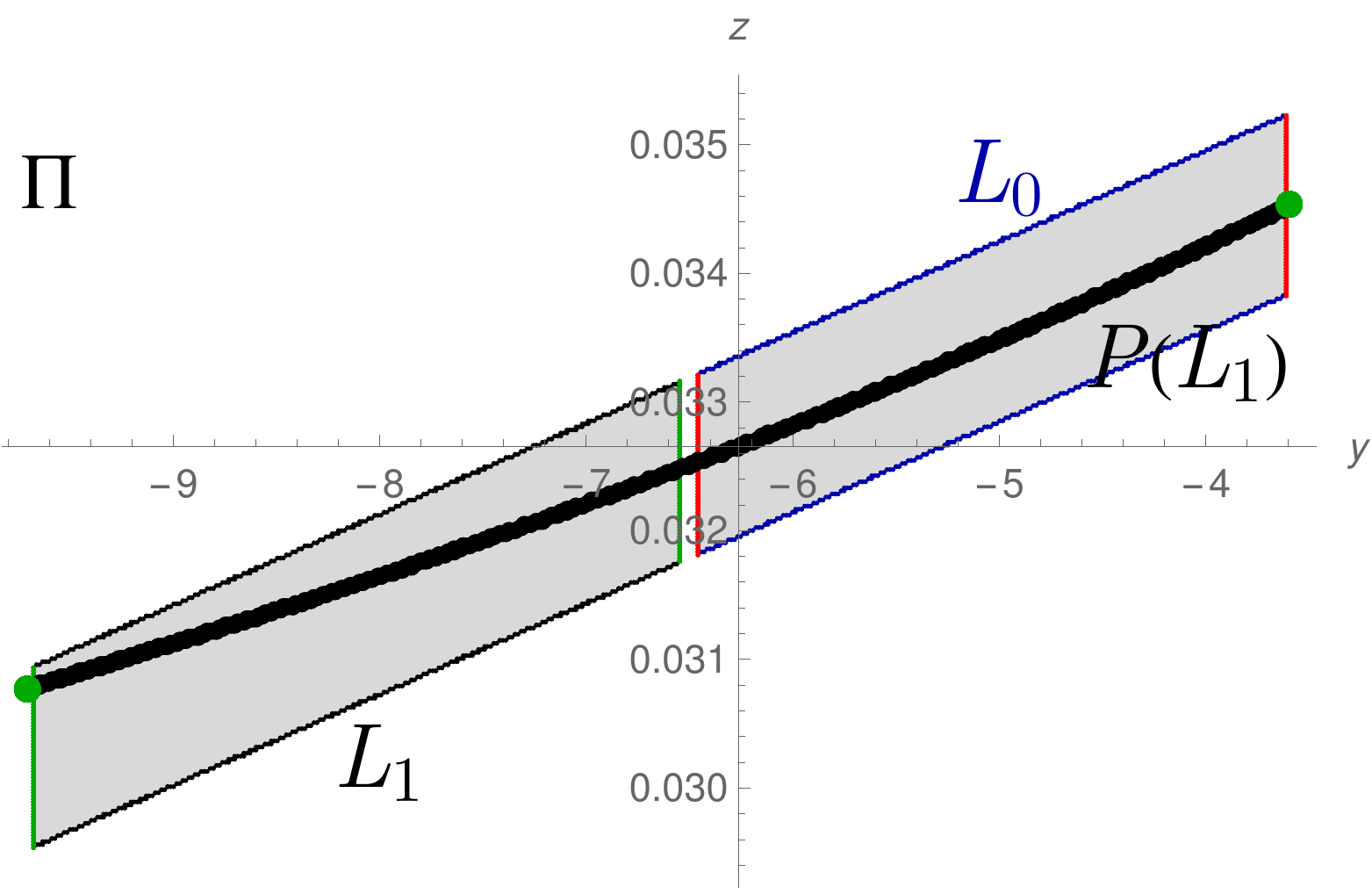}
	\caption{\label{fig:3-per}The sets $L_0$, $L_1$ and their images through $P$. The axes' origin is moved to the point $p_2^3$. The sets fulfil $L_0\overset{P}{\Longrightarrow} L_1 \overset{P}{\Longrightarrow} L_1
\overset{P}{\Longrightarrow} L_0$ in an affine coordinate system.}
\end{figure}

Denote $P_C =C^{-1}\circ P\circ C$. Then the following horizontal covering relations occur (from computer-assisted proof \cite{proof}, Case 2 of the program \texttt{01-Roessler\_a525.cpp}, see also the outline in Appendix, Subsec.\ \ref{ss:prog_covering}):
\[
N_0\overset{\text{\raisebox{2pt}{$P_C$}}}{\Longrightarrow} N_1 \overset{\text{\raisebox{2pt}{$P_C$}}}{\Longrightarrow} N_1
\overset{\text{\raisebox{2pt}{$P_C$}}}{\Longrightarrow} N_0.
\]
Now, all periods for $P_C$ can be obtained as the loops of horizontal $P_C$-covering, as in Theorem \ref{th:periodic}:
\begin{itemize}
	\item for $n=1$: the existence of a stationary point follows from the self-covering
	\[N_1 \overset{P_C}{\Longrightarrow} N_1\text{;}\]
	\item for $n\geq 2$: the existence of an $n $-periodic orbit follows from the chain of covering relations:
	\[
	N_0 \overset{P_C}{\Longrightarrow} N_1\underbrace{\overset{P_C}{\Longrightarrow} N_1 \dots \overset{P_C}{\Longrightarrow}N_1}_{\text{`}\overset{P_C}{\Longrightarrow} N_1\text{' } n-2 \text{ times}}  \overset{P_C}{\Longrightarrow} N_0.
	\]
Note that in this case $n$ must be the fundamental period of the point.
\end{itemize}
Finally, observe that an $n $-periodic orbit for $P_C$ defines an $n $-periodic orbit for $P$.
\end{proof}

\section{Periodic orbits for $a=4.7$}

Consider now the system \eqref{eq:rossler} with $a=4.7$, $b=0.2$, that is
\begin{equation}\label{eq:rossler47}
\begin{cases}
x'=-y-z,
\\
y'=0.2\, y+x,
\\
z'=z (x-4.7)+0.2
\end{cases}
\end{equation}

Our goal is to establish the following result.

\begin{figure}[h]
	\includegraphics[height=6cm]{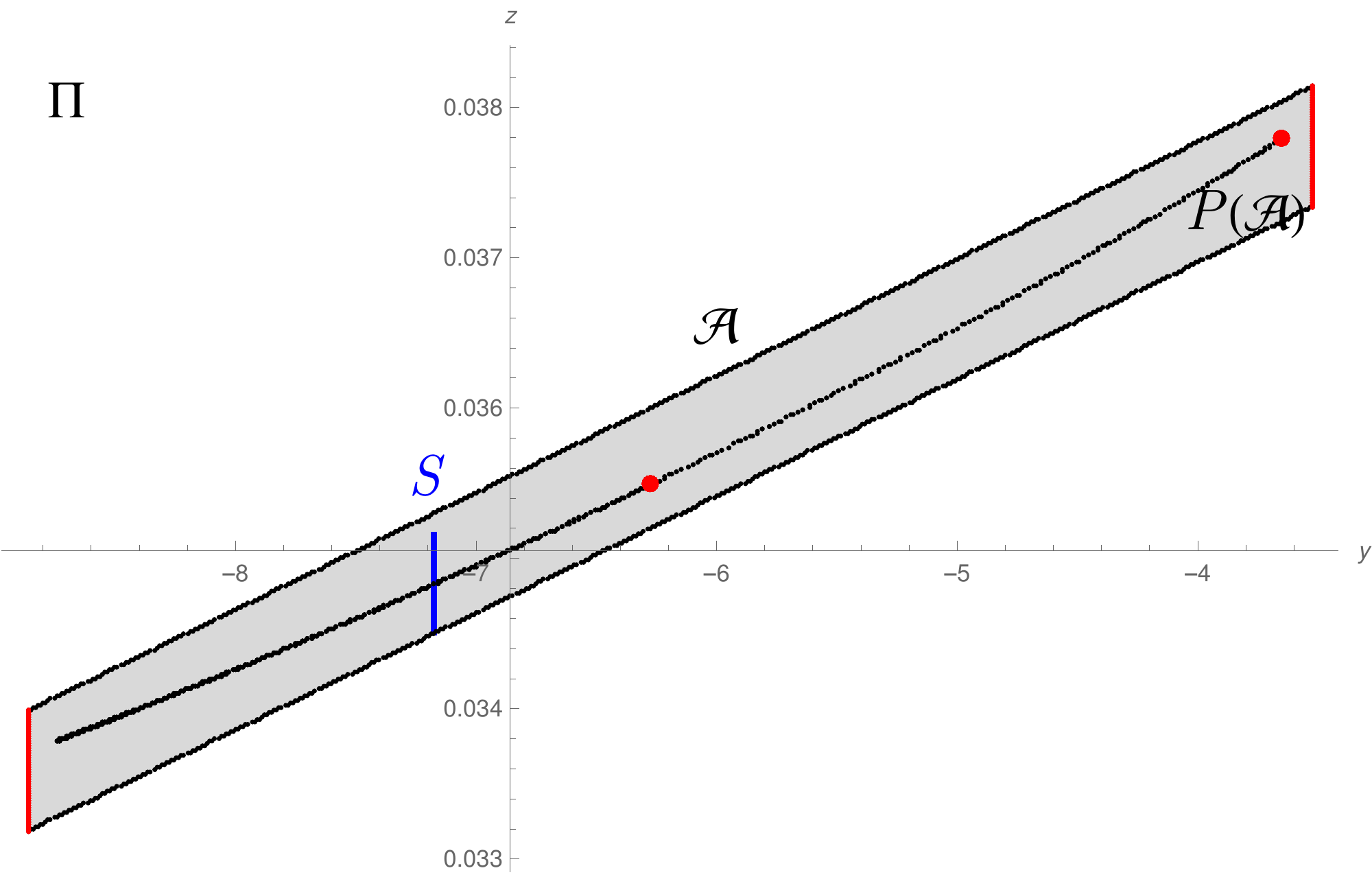}
	\caption{\label{fig:5-attr_container}The forward-invariant set $\mathcal{A}$ containing an attractor for the system \eqref{eq:rossler47}, and its image through $P$. The blue rectangle $S$ contains a stationary point.}
\end{figure}

\begin{theorem}\label{th:r5}
Consider a parallelogram $\mathcal{A}$ in the $(y,z)$ coordinates on the section $\Pi$ (see Fig. \ref{fig:5-attr_container}):
	\begin{equation}\label{eq:A}
		\mathcal{A} = \begin{bmatrix}-6.19384 \\ 0.0356629\end{bmatrix}
		+	
		\begin{bmatrix}-1. & 0.000777754 \\ -0.000777754 & -1.\end{bmatrix}
		\cdot
		\begin{bmatrix}\pm 2.66856 \\ \pm 4\cdot 10^{-4}\end{bmatrix}.	
	\end{equation}

	Then $\mathcal{A}$ is forward-invariant for the map $P$, that is
	\begin{equation}
	P(\mathcal{A})\subset \mathcal{A}\text{,}  \label{eq:PA}
	\end{equation}
and the R\"ossler system \eqref{eq:rossler47} has $n$-periodic orbits for any $n \in \mathbb{N}\setminus \{3\}$, passing through $\mathcal{A}$, and it does not have any $3$-periodic orbit there.
\end{theorem}
Let us comment first the statement about the existence of the forward-invariant $\mathcal{A}$. Such a set was not mentioned in the statement of Theorem~\ref{th:r3} because we had proved there the existence of all periods. Now beside the existence of some periodic points we also want to exclude the period $3$, and for this we need to be precise about where this exclusion happens.

The proof relies on several lemmas.

\begin{figure}[h]
	\includegraphics[height=9cm]{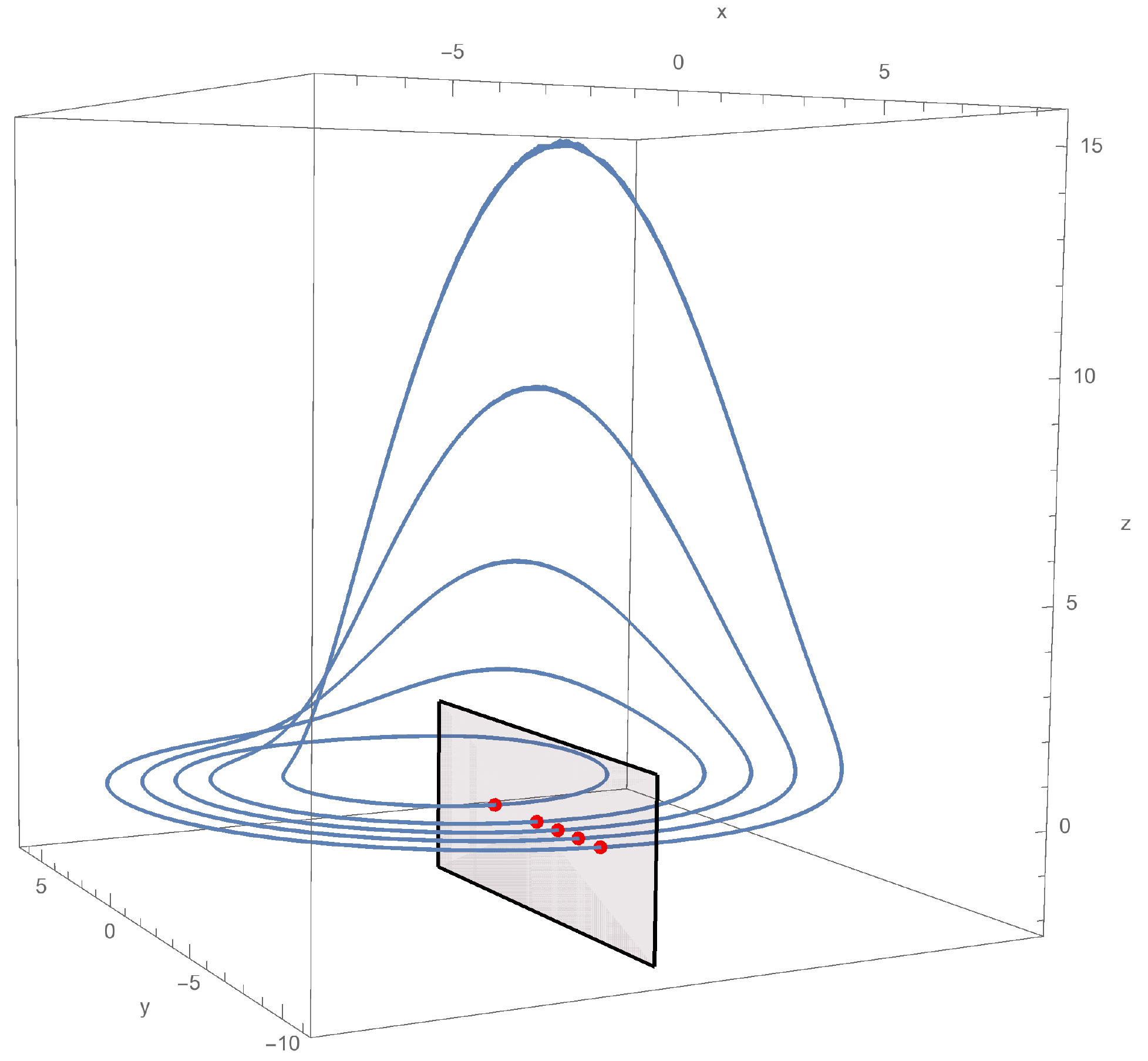}
	\caption{\label{fig:attr5}The attracting $5$-periodic orbit for the system \eqref{eq:rossler47}.}
\end{figure}

From the bifurcation diagram (Fig. \ref{fig:bif}) it is apparent that there exists a $5$-periodic orbit for the system \eqref{eq:rossler47} (see Fig. \ref{fig:attr5}), contained in $\mathcal{A}$. As previously, one proves it by interval Newton method \cite{proof} (Case 1 of the program \texttt{02-Roessler\_a47.cpp}, see also Appendix, Subsec.\ \ref{ss:prog_periodic}).

\begin{lemma}\label{lem:5per}
The Poincar\'e map $P$ of the system \eqref{eq:rossler47} on the section $\Pi$ has a 5-periodic orbit $\mathcal{O}^5=\{p_1^5,p_2^5,p_3^5,p_4^5,p_5^5\}$, contained in the following rectangles in the $(y,z)$ coordinates on the section $\Pi$:
\begin{equation}
\begin{aligned}\label{eq:5-per}
		p_1^5 \in & -3.885277116_{910041}^{888829} \times 0.037558391444_{32487}^{85094}, \\	
		p_2^5 \in & -6.8582604471_{62484}^{26429} \times 0.03505366666_{495363}^{561609}, \\
		p_3^5 \in & -7.7666312453_{92379}^{48371} \times 0.034417133928_{18681}^{63033}, \\	
		p_4^5 \in & -5.895584354_{611225}^{509201} \times 0.03578591178_{706747}^{835873}, \\
		p_5^5 \in & -8.7223960200_{73123}^{49997} \times 0.033796299364_{04972}^{3551}.
		\end{aligned}
\end{equation}
\end{lemma}

\begin{lemma}\label{lem:r5-432}
The R\"ossler system \eqref{eq:rossler47} has $n$-periodic orbits for any $n \in \mathbb{N}\setminus \{4,3,2\}$.
\end{lemma}

\noindent\textit{Heuresis.}

\begin{figure}[h]
	\includegraphics[width=0.45\textwidth]{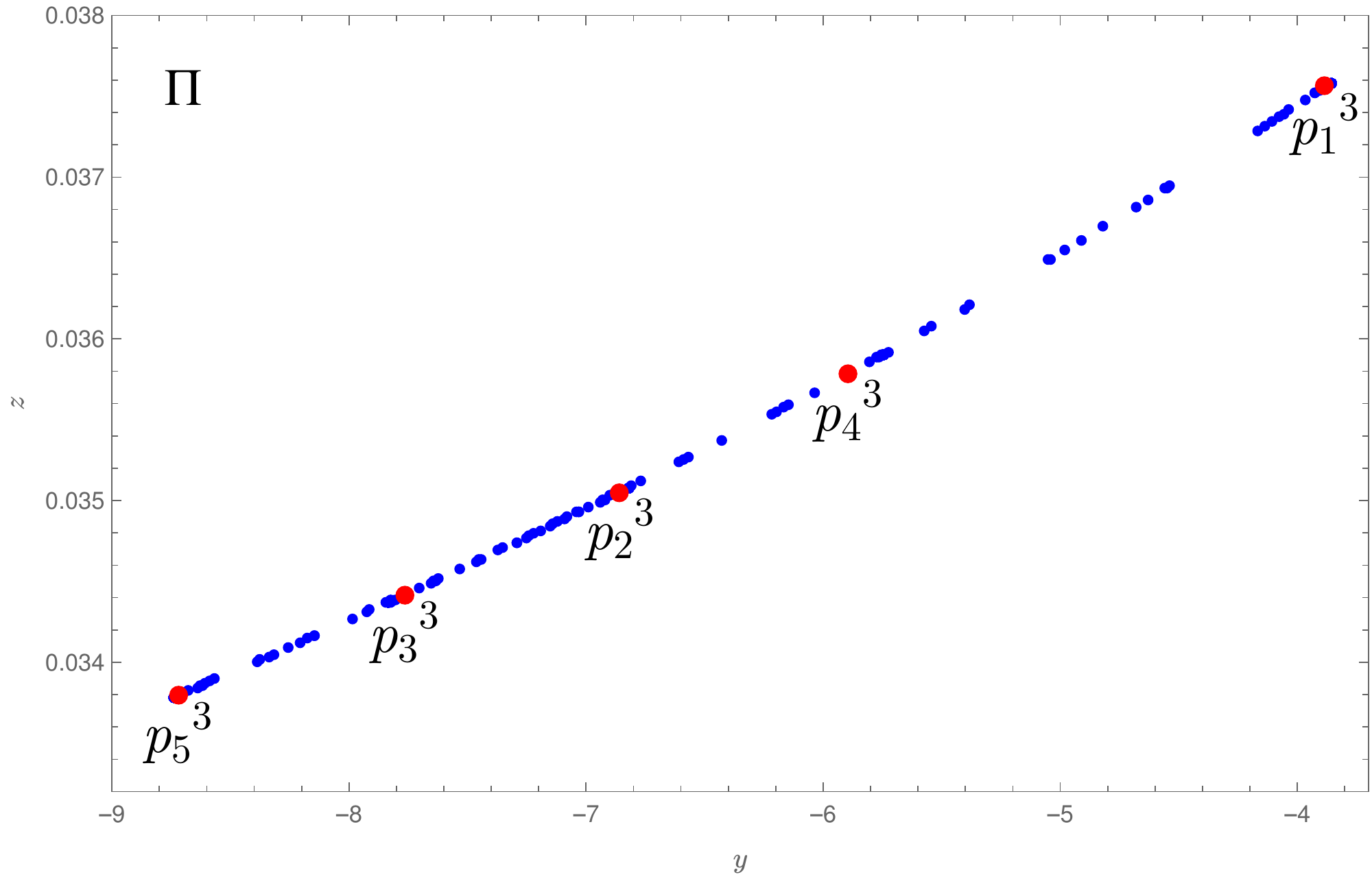}	
	\includegraphics[width=0.45\textwidth]{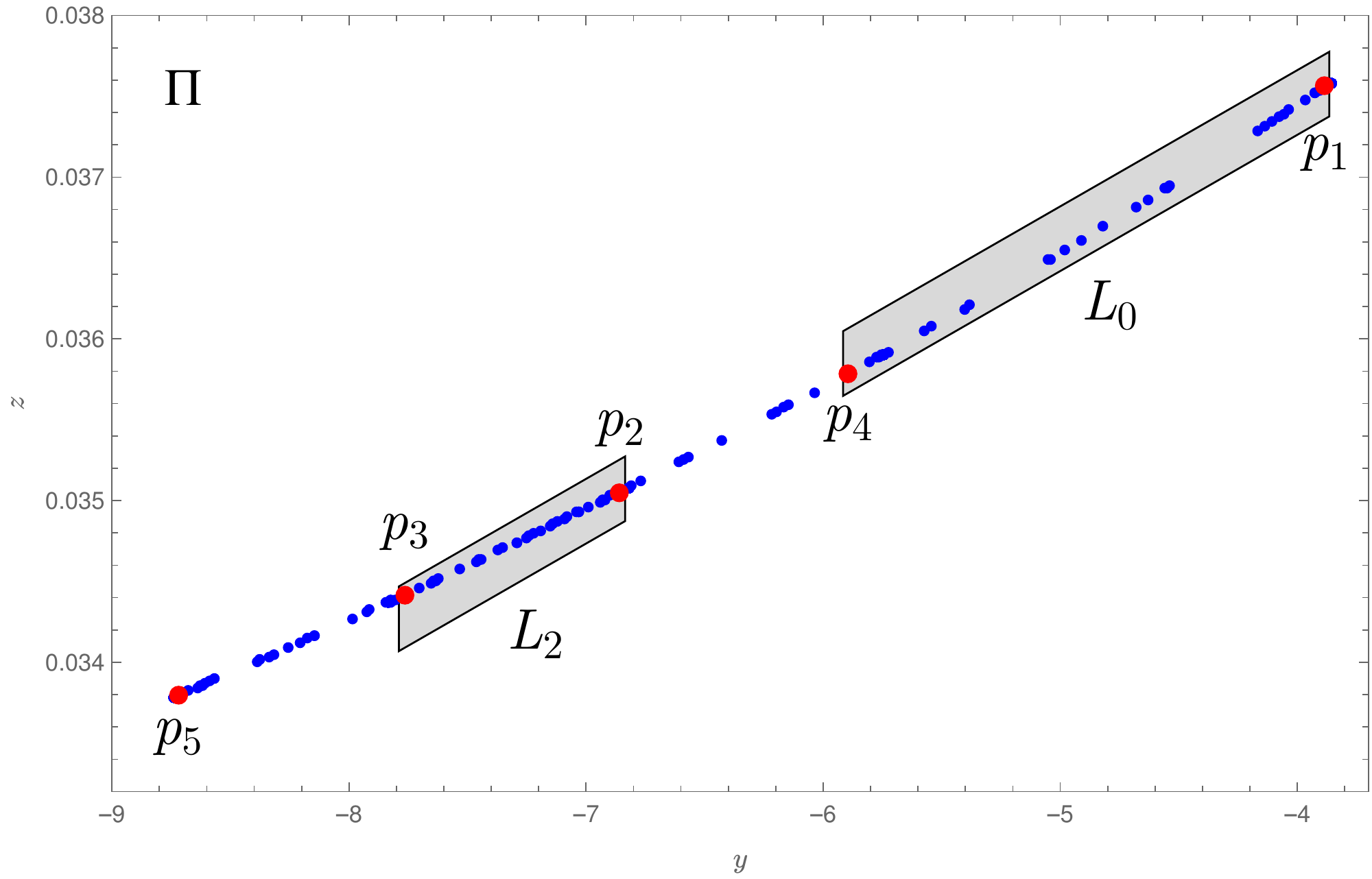}
	\caption{\label{fig:intervalp12345}\underline{To the left:} The fragment of $\Pi$ section containing the $5$-periodic orbit $\mathcal{O}^5$, marked in red. Some orbits attracted by  $\mathcal{O}^5$ allow to estimate the shape of $\mathcal{M}$.
	\newline
	\label{fig:L_vs_p5}\underline{To the right:} The location of the sets $L_0$, $L_2$ relative to the orbit $\mathcal{O}^5$.}
\end{figure}

Similarly as in the proof of Theorem \ref{th:r3}, consider a hypothetical 1-dimensional manifold $\mathcal{M}\subset \Pi$ containing $\mathcal{O}^5$ and its self-map, for which  the Poincar\'e map $P$ is a 2-dimensional perturbation (see Fig. \ref{fig:intervalp12345}).

\begin{figure}[h]
	\includegraphics[width=0.5\textwidth]{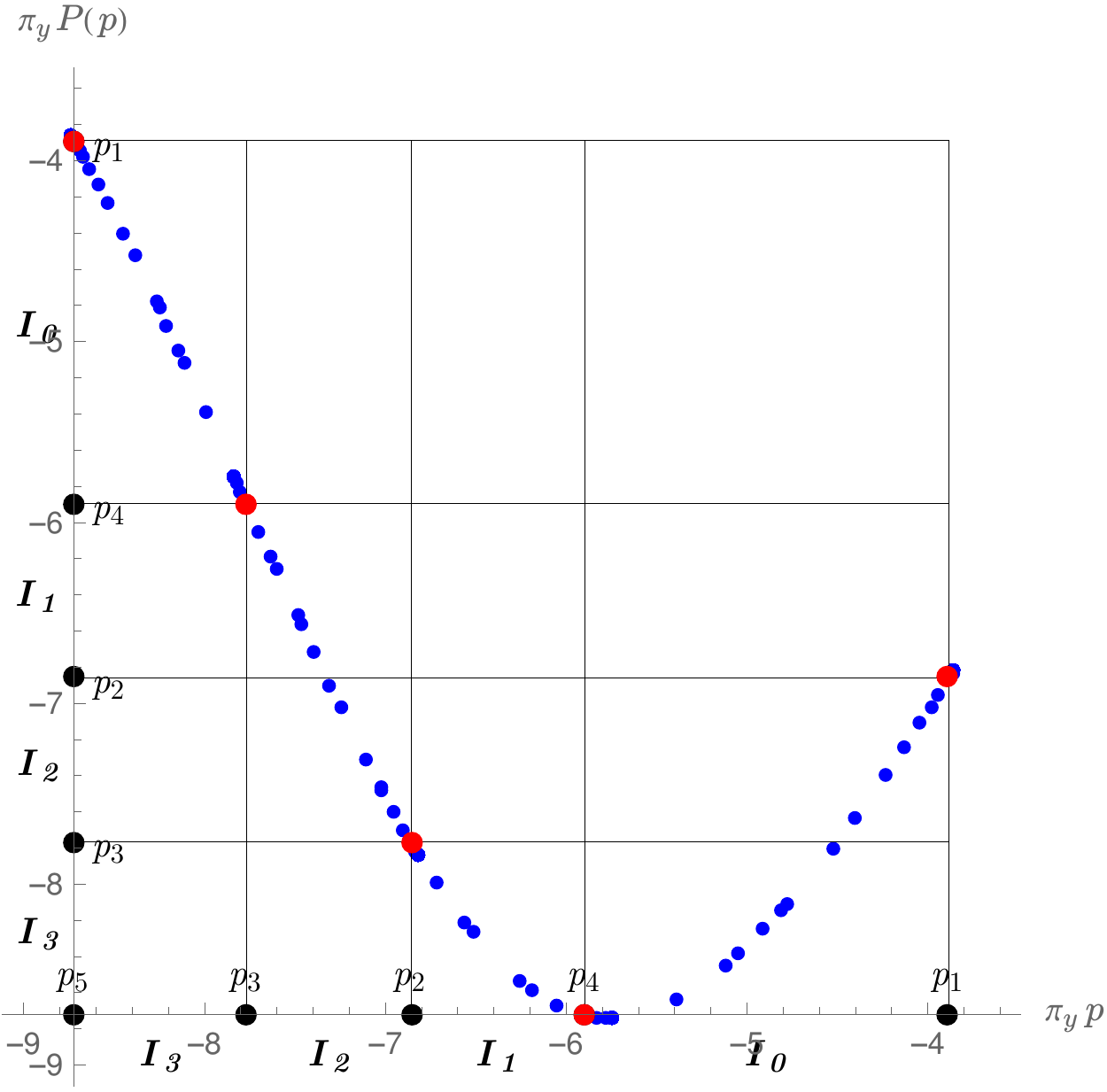}
	\caption{\label{fig:py_vs_Ppy5}The images of the points near $\mathcal{M}$ through the model map $\mathcal{P}$.}
\end{figure}

Let us now parameterize $\mathcal{M}$ by the $y$ coordinate and consider, as previously, the `plot' of the model map  $\mathcal{P}: \mathcal{M} \to \mathcal{M}$.  Denote the four `segments' of $\mathcal{M}$ by $I_0$, $I_1$, $I_2$ and $I_3$, as depicted on Fig. \ref{fig:py_vs_Ppy5}.
In this case, $\mathcal{M}$'s segments fulfil the following diagram of one-dimensional covering relations:
\begin{equation}\label{eq:diag}
\xymatrix{
I_0
	\ar@<.5ex>[d]^{\mathcal{P}}
	\ar[r]^{\mathcal{P}}
& I_2 \ar[d]^{\mathcal{P}} \ar@(ru,rd)^{\mathcal{P}}
\\
I_3 \ar@<.5ex>[u]^{\mathcal{P}}
& I_1 \ar[l]_{\mathcal{P}}
}
\end{equation}
In particular, the following chain of covering relations occurs:
\[
I_0 \overset{\mathcal{P}}{\longrightarrow}I_2 \overset{\mathcal{P}}{\longrightarrow}I_2 \overset{\mathcal{P}^3}{\longrightarrow} I_0.
\]
Similarly as before, we find two-dimensional sets $L_0$, $L_2$ on $\Pi$, close to the segments $I_0$, $I_2$, which are expected to fulfil the analogous horizontal covering relations. One can compare their positions to the orbit $\mathcal{O}^5$ on Fig.\ \ref{fig:L_vs_p5}, to the right.
Further we will see the proof of these horizontal covering relations.

\begin{proof}~

Let  $M = \begin{bmatrix}
		-1 & 0.000842495 \\-0.000842495 & -1
	\end{bmatrix}$ and $p_2^5=\begin{bmatrix}
	-6.858260447127058 \\ 0.03505366666527084
	\end{bmatrix}$. Define an affine map $C$ on $\mathbb{R}^2$ by $	C(x)=Mx+p_2^5$ and denote  $P_C =C^{-1}\circ P\circ C$. As before, the map $C$ is chosen to straighten approximately the set $\mathcal{M}$.
Consider now two h-sets $N_0,N_2\subset \mathbb{R}^2$:

\begin{equation}
\begin{aligned}
	N_0 &= [-1.96783 \pm 1.02]\times [\pm 2\cdot 10^{-4}]
	\text{,}
	\\
	N_2 &= [	0.454186	\pm 0.476895]\times [\pm 2\cdot 10^{-4}].
\end{aligned}
\end{equation}

Denote their images through $C$ by $L_0 = C(N_0)$ and $L_2 = C(N_2)$, and consider these parallelograms as sets on the section $\Pi$ (see Fig. \ref{fig:L_vs_p5}, to the right). Compare also to Fig.\ \ref{fig:5-per}, where the images of $L_0$, $L_2$ on the $\Pi$ section are depicted.
\begin{figure}[h]
	\includegraphics[width=0.48\textwidth]{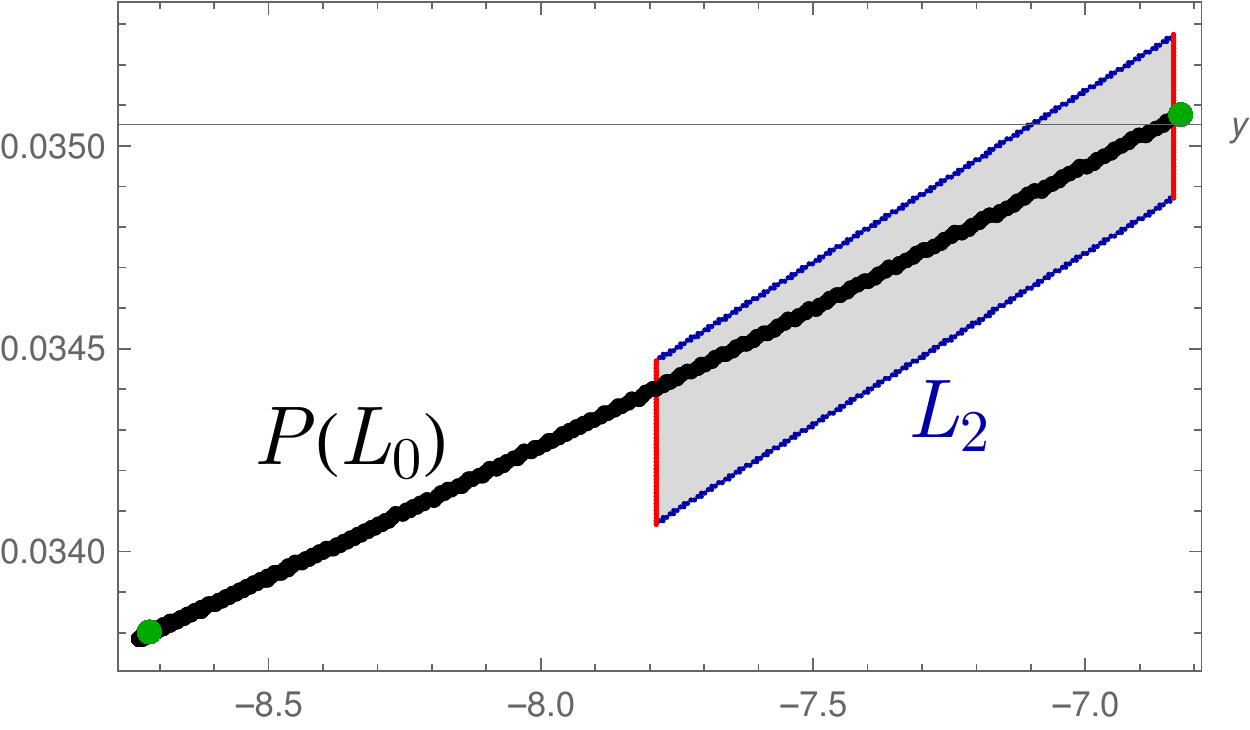}	
	\includegraphics[width=0.48\textwidth]{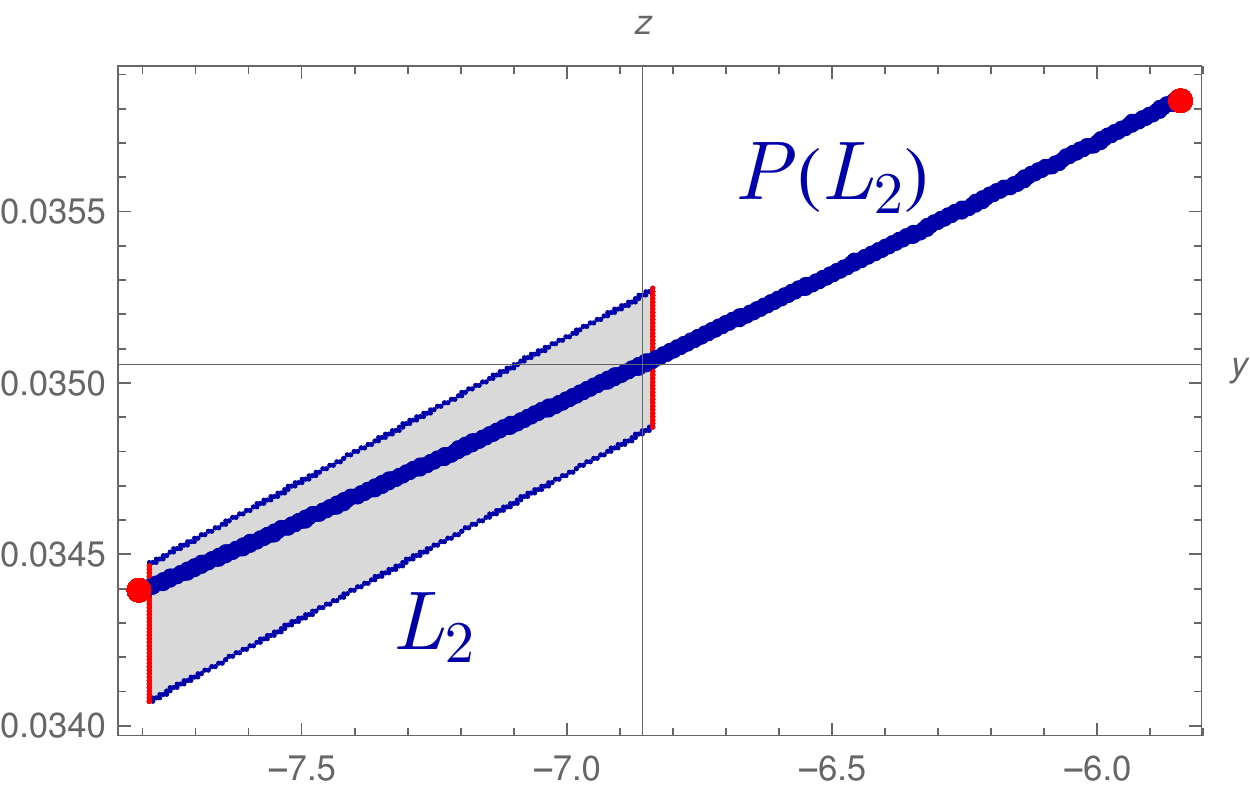}
	\includegraphics[width=0.48\textwidth]{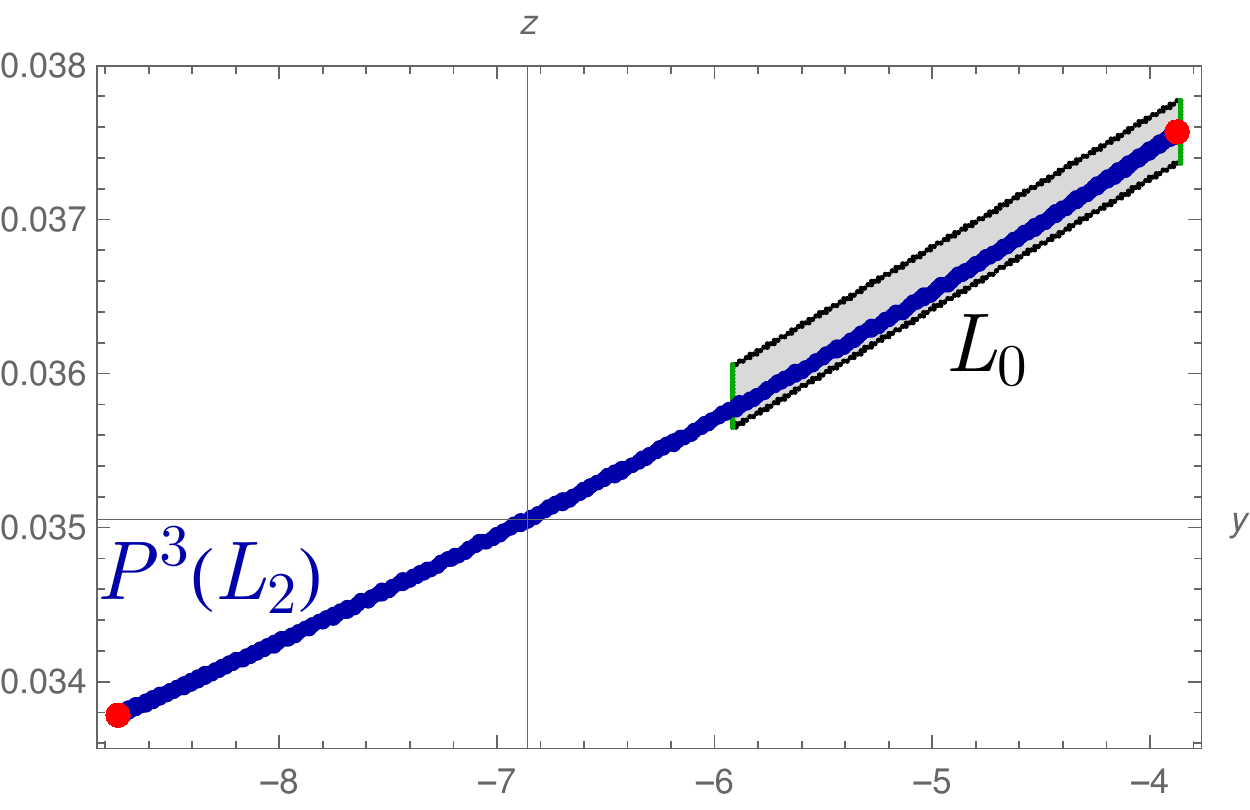}
	\caption{\label{fig:5-per}The sets $L_0$, $L_2$ and their images through $P$ and $P^3$. The axes' origin is moved to the point $p_2^5$. The sets fulfil $L_0\overset{P}{\Longrightarrow} L_2 \overset{P}{\Longrightarrow} L_2
\overset{P^3}{\Longrightarrow} L_0$ in an affine coordinate system.}
\end{figure}

Let $P_C =C^{-1}\circ P\circ C$. Then the following horizontal covering relations occur (from computer-assisted proof \cite{proof}, Case 2 of the program \texttt{02-Roessler\_a47.cpp}, see also the outline in Appendix, Subsec.\ \ref{ss:prog_covering}):
\[
N_0\overset{\text{\raisebox{2pt}{$P_C$}}}{\Longrightarrow} N_2 \overset{\text{\raisebox{2pt}{$P_C$}}}{\Longrightarrow} N_2
\overset{\text{\raisebox{2pt}{$P_C^3$}}}{\Longrightarrow} N_0.
\]
Now, we obtain all desired periods for $P_C$ as the loops of horizontal $P_C$-covering, as in Theorem \ref{th:periodic}:
\begin{itemize}
	\item for $n=1$: the existence of a stationary point follows from the self-covering
	\[N_2 \overset{P_C}{\Longrightarrow} N_2\text{;}\]
	\item for $n\geq 5$: the existence of an $n$-periodic orbit follows from the chain of covering relations:
	\[
	N_0 \overset{P_C}{\Longrightarrow} N_2\underbrace{\overset{P_C}{\Longrightarrow} N_2 \dots \overset{P_C}{\Longrightarrow}N_2}_{\text{`}\overset{P_C}{\Longrightarrow} N_2\text{' } n-4 \text{ times}}  \overset{P_C^3}{\Longrightarrow} N_0.
	\]
Also in this case $n$ must be the fundamental period of the point.
\end{itemize}
An $n$-periodic orbit for $P_C$ defines an $n$-periodic orbit for $P$, so we get all periods $n$ for $n\in \mathbb{N}\setminus\{2,3,4\}$.
\end{proof}

To complete the study, we find also $2$- and $4$-periodic orbits. Note that from Fig.\ \ref{fig:py_vs_Ppy5} and diagram \eqref{eq:diag} we expect that for the model map $\mathcal{P}$:
\[
	I_0 \overset{\mathcal{P}^2}{\longrightarrow} I_0
	\text{, \quad and \quad }
	I_1 \overset{\mathcal{P}^4}{\longrightarrow} I_1
	\text{,}
\]
therefore we may look for a $2$-periodic point on the segment $I_0$ and for a $4$-periodic point on $I_1$.

\begin{lemma}\label{lem:r5+24}~

	\begin{enumerate}
	\item The Poincar\'e map $P$ of the system \eqref{eq:rossler47} on the section $\Pi$ has a 2-periodic orbit $\mathcal{O}^2=\{p_1^2,p_2^2\}$, contained in the following rectangles in the $(y,z)$ coordinates on the section $\Pi$:
		\begin{equation}
		\begin{aligned}\label{eq:2-per}
				p_1^2 \in & -4.88392425874_{3846}^{2264} 				\times 0.0366312810972_{0363}^{9599}, \\	
				p_2^2 \in & -8.22095155223_{5825}^{3453} 				\times 0.0341162008426_{8562}^{9375}.
		\end{aligned}
		\end{equation}
	\item The Poincar\'e map $P$ of the system \eqref{eq:rossler47} on the section $\Pi$ has a 4-periodic orbit $\mathcal{O}^4=\{p_1^4,p_2^4,p_3^4,p_4^4\}$, contained in the following rectangles in the $(y,z)$ coordinates on the section $\Pi$:
		\begin{equation}
		\begin{aligned}\label{eq:4-per}
				p_1^4 \in & -6.3322511801_{91433}^{86209} 				\times 0.035445799547_{48502}^{86024}, \\	
				p_2^4 \in & -8.4703436541_{25783}^{19862} 				\times 0.033955548564_{30393}^{40434}, \\
				p_3^4 \in & -4.3626672599_{1514}^{017} 				\times 0.037102455850_{53683}^{75971}, \\	
				p_4^4 \in & -7.5716695403_{62099}^{41765} 				\times 0.0345497016342_{597}^{7861}.
		\end{aligned}
		\end{equation}
	\end{enumerate}
\end{lemma}

\begin{proof}
Computer-assisted \cite{proof}, Cases 3 and 4 of the program \texttt{02-Roessler\_a47.cpp}. See also the outline in Appendix, Subsec.\ \ref{ss:prog_periodic}.
\end{proof}

\noindent\textit{Proof of Theorem~\ref{th:r5}.}

First we verify the forward invariance of $\mathcal{A}$: \cite{proof}, Case 5 of the program \texttt{02-Roessler\_a47.cpp} (see also the outline in Appendix, Subsec.\ \ref{ss:prog_non3}).

	From Lemmata \ref{lem:r5-432} and \ref{lem:r5+24} we have already $n$-periodic orbits for $n\in \mathbb{N}\setminus\{3\}$ in $\mathcal{A}$. It remains to  prove that there is no $3$-periodic orbit for $P$ in $\mathcal{A}$. We verify that by the interval Newton method applied to $P$ and $P^3$: \cite{proof}, Case 6 of the program \texttt{02-Roessler\_a47.cpp} (see also the outline in Appendix, Subsec.\ \ref{ss:prog_non3}).
\qed

\section{Appendix. Outlines of computer-assisted proofs \cite{proof}}\label{app}

We use the CAPD library for C++ \cite{capd-article} containing, in particular, modules for interval arithmetic, linear algebra, differentiation and integration of ODEs. The important utility of the library is calculating Poincar\'e maps rigorously: we are able to enclose the image of an interval vector (a rectangle $R$, in our case) through $P$ in an interval vector, denoted below by $[P(R)]$. Applying suitable affine coordinate systems helps to reduce the wrapping effect.
If that is not enough and the estimated images are too large, we also often divide the sets into grids of smaller boxes and use the fact that the image of the whole set must be contained in the sum of small boxes' images.

The programs have been tested under Linux Mint 18.1 with \texttt{gcc} 
compiler. They use the CAPD library ver. 5.0.6.
All cases execute within seconds on a laptop type computer with Intel Core i7 2.7GHz $\times$ 2 processor
(Case 6 of the program \texttt{02-Roessler\_a47.cpp} takes the longest time: no more than 10 seconds).

\subsection{Detecting stationary points for a Poincar\'e map's $n$th iterate $P^n$}\label{ss:prog_periodic}
The outline refers to the following proofs:
\begin{itemize}
	\item existence of a $3$-periodic orbit for $P$, $a=5.25$ (Lemma \ref{lem:3-per})\\
	-- Case 1 of the program \texttt{01-Roessler\_a525.cpp};
	\item existence of a $5$-periodic orbit for $P$, $a=4.7$ (Lemma \ref{lem:5per})\\
	-- Case 1 of the program \texttt{02-Roessler\_a47.cpp};
	\item existence of $2$- and $4$-periodic orbits for $P$, $a=4.7$ (Lemma \ref{lem:r5+24})\\
	-- Cases 3 and 4 of the program \texttt{02-Roessler\_a47.cpp}.
\end{itemize}

\textit{Outline of the computer-assisted proof}.
\begin{enumerate}
	\item We define the system and a starting point $x_0$ close to the expected periodic point;
	\item We apply the function \texttt{anyStationaryPoint}, which:
	\begin{itemize}
		\item initially iterates the interval Newton operator (INO) 
		\[
		N(x_0,X)=x_0 - [DP^n(X)]^{-1} P^n(x_0)
		\]
		in a small neighbourhood $X$ of $x_0$ to get a better starting point (this step is not necessary, but we finally obtain a better approximation of the periodic point);
		\item finally applies INO on the small neighbourhood of the better starting point: if the image $N$ is contained in the interior of the neighbourhood, then a unique stationary point exists in it.
		\item returns $N$.
	\end{itemize}
	\item Finally we apply a given iteration of $P$ to $N$, to print the whole orbit.
\end{enumerate}

\subsection{Horizontal covering relations between h-sets on the section $\Pi$ for the Poincar\'e map $P$ or its higher iterate}\label{ss:prog_covering}~

The outline refers to the following proofs:
\begin{itemize}
	\item the chain of covering relations $N_0\overset{\text{\raisebox{2pt}{$P_C$}}}{\Longrightarrow} N_1 \overset{\text{\raisebox{2pt}{$P_C$}}}{\Longrightarrow} N_1
\overset{\text{\raisebox{2pt}{$P_C$}}}{\Longrightarrow} N_0$ (Theorem \ref{th:r3}) \\
-- Case 2 of the program \texttt{01-Roessler\_a525.cpp};
	\item the chain of covering relations $N_0\overset{\text{\raisebox{2pt}{$P_C$}}}{\Longrightarrow} N_2 \overset{\text{\raisebox{2pt}{$P_C$}}}{\Longrightarrow} N_2
\overset{\text{\raisebox{2pt}{$P_C^3$}}}{\Longrightarrow} N_0$ (Theorem \ref{th:r5}) \\
-- Case 2 of the program \texttt{02-Roessler\_a47.cpp};
\end{itemize}

\textit{Outline of the computer-assisted proof}.
\begin{enumerate}
	\item We define the system and the h-sets with their affine coordinate systems;
	\item We apply the function \texttt{covers2D}, which checks the conditions from Definition \ref{def:cov} for the estimated images of the h-sets and their vertical edges;
	\item If the conditions are fulfilled, then the function returns \texttt{true}.
\end{enumerate}
See also a similar method for 2-dimensional covering  described in \cite{AGPZ}.

\subsection{Forward invariance of the set $\mathcal{A}$}\label{ss:prog_A}~

The outline refers to the case of the Poincar\'e map $P$, $a=4.7$ and the set $\mathcal{A}$ from Theorem \ref{th:r5}) -- Case 5 of the program \texttt{02-Roessler\_a47.cpp}.

\textit{Outline of the computer-assisted proof}.
\begin{enumerate}
	\item We define the set $\mathcal{A}$ with its affine coordinate system;
	\item We apply the function \texttt{inside}, which divides $\mathcal{A}$ into small boxes $A_i$ and checks if each $A_i$ is mapped to the interior of  $\mathcal{A}$.
	\item If the above condition is fulfilled, then the function returns \texttt{true}.
\end{enumerate}

\subsection{Non-existence of the $3$-periodic point for the Poincar\'e map $P$ in the set $\mathcal{A}$}\label{ss:prog_non3}~

The outline refers to the case of the Poincar\'e map $P$, $a=4.7$ and the set $\mathcal{A}$ from Theorem \ref{th:r5}) -- Case 6 of the  program \texttt{02-Roessler\_a47.cpp}.

\textit{Outline of the computer-assisted proof}.
	\begin{enumerate}
	\item We define the set $\mathcal{A}$ with its affine coordinate system;
	\item We apply the function \texttt{whatIsNotMappedOutside}, which:
	\begin{itemize}
		\item divides $\mathcal{A}$ into $500\times 10$ small boxes $A_i$
		\item checks if each image $[P^3(A_i)]\cap A_i = \varnothing$
		\item returns a set $S$, which is an interval closure of the sum of all $A_i$'s such that $[P^3(A_i)]\cap A_i \neq \varnothing$ (see on Fig.\ \ref{fig:5-attr_container}):
			\[
			S = -7.1_{86544568881281}^{65195528898398}\times
			0.03_{44908202443027}^{522742411001666}.
			\]
	\end{itemize}	
		\item We check by interval Newton method (function \texttt{newtonDivided}) that $S$ contains a unique stationary point for $P^3$; it still may be in fact a stationary point for $P$.
		\item We check by \texttt{newtonDivided} that $S$ contains a unique stationary point also for $P$ and therefore $S$ does not contain any points of fundamental period $3$.
	\end{enumerate}

\bibliographystyle{plain}
\bibliography{roessler_bib}

\end{document}